\documentclass[opre,nonblindrev]{informs3ArXiv} 

\OneAndAHalfSpacedXI 

\usepackage{endnotes}
\let\footnote=\endnote

\newtheorem{fact}{Fact}

\usepackage{natbib}
 \bibpunct[, ]{(}{)}{,}{a}{}{,}%
 \def\bibfont{\small}%
 \def\bibsep{\smallskipamount}%
 \def\bibhang{24pt}%
\usepackage{rotating}
\usepackage{algorithm}
\usepackage{algorithmic}
\usepackage{amsbsy}
\usepackage{graphicx}
\usepackage{caption}

\TheoremsNumberedThrough     
\ECRepeatTheorems

\EquationsNumberedThrough    

\begin{document}


\RUNAUTHOR{Thoppe, Borkar and Garg}

\RUNTITLE{Quadratic Optimization in High Dimensions}

\TITLE{Greedy Block Coordinate Descent (GBCD) Method for High Dimensional Quadratic Programs}

\ARTICLEAUTHORS{%
\AUTHOR{Gugan Thoppe}
\AFF{School of Technology and Computer Sc., Tata Institute of Fundamental Research, Mumbai 400005, India, \EMAIL{gugan@tcs.tifr.res.in}} 
\AUTHOR{Vivek S.\ Borkar}
\AFF{Department of Electrical Engineering, Indian Institute of Technology Bombay, Mumbai 400076, India, \EMAIL{borkar.vs@gmail.com}}
\AUTHOR{Dinesh Garg}
\AFF{IBM Research - India, Manyata Embassy Business Park, Bangalore 560045, India, \EMAIL{garg.dinesh@in.ibm.com}}

} 

\ABSTRACT{%
High dimensional unconstrained quadratic programs (UQPs) involving massive datasets are now common in application areas such as web, social networks, etc. Unless computational resources that match up to these datasets are available, solving such problems using classical UQP methods is very difficult. This paper discusses alternatives. We first define high dimensional compliant (HDC) methods for UQPs---methods that can solve high dimensional UQPs by adapting to available computational resources. We then show that the class of block Kaczmarz and block coordinate descent (BCD) are the only existing methods that can be made HDC. As a possible answer to the question of the `best' amongst BCD methods for UQP, we propose a novel greedy BCD (GBCD) method with serial, parallel and distributed variants. Convergence rates and numerical tests confirm that the GBCD is indeed an effective method to solve high dimensional UQPs. In fact, it sometimes beats even the conjugate gradient.
}%



\maketitle

%

\section{Introduction} \label{sec:Intro}
Machine learning and statistics problems arising in application areas such as the web, social networks, e-commerce, m-commerce, etc. are often extremely large in size. As described in \cite{boyd2011distributed}, these problems typically share two characteristics---(1) the input datasets are staggeringly large, consisting of millions or billions of training examples, and (2) the data describing each example is itself high-dimensional. While many of these problems are still convex optimization programs, the above two {\em high dimensional} characteristics bring in new challenges concerning management of data. The worst is when the available computational setup to solve these problems is itself minimal and simplistic. As a result, methods which can solve large scale optimization problems by adapting to the resources available are now of central importance. We will refer to such methods as {\em high dimensional compliant} (HDC) methods. In this paper, we aim to build HDC solution methods for an unconstrained quadratic program (UQP)---the simplest convex optimization problem.

Recent literature, see \cite{nesterov2012efficiency}, \cite{richtarik2012parallel}, \cite{needell2014paved}, etc., suggests that block Kaczmarz (BK) and block coordinate descent (BCD) class of methods may be the only ones suited to tackle large scale optimization. A rigorous verification of this belief has, however, been missing. In this paper, we first define the notion of HDC solution methods for UQP. We then show that BK and BCD are the only existing UQP solution methods which can be made HDC. In fact, we prove that even natural generalizations of the BK and BCD philosophies are not HDC. This then brings us to the question of which is the `best' amongst BK and the `best' amongst BCD methods for a UQP. As a possible answer to the latter, we propose a novel deterministic, but adaptive, greedy BCD (GBCD) method with serial, parallel, and distributed variants. Convergence rates and numerical tests confirm that GBCD is indeed an effective method to solve high dimensional UQPs. In fact, it sometimes beats even the conjugate gradient.

\subsection{Preliminaries}
Throughout this paper, a UQP will essentially be assumed to have the general form
\begin{equation}
\label{eqn:UQP}
\underset{x \in \mathbb{R}^n}{\min} \hspace{1ex} f(x) = \frac{1}{2}x^tPx - x^tq + r,
\end{equation}
where $P \in \mathbb{S}^{n}_{++} \subset \mathbb{R}^{n \times n},$ i.e., the matrix $P$ is symmetric and positive definite, $q \in \mathbb{R}^n$ and $r \in \mathbb{R}.$  We will often refer to $P$ as the input matrix. Since the gradient of the quadratic function $f$  is
\begin{equation}
\label{eqn:gradUQP}
\nabla f(x) = Px - q
\end{equation}
and $P \in \mathbb{S}^{n}_{++},$ it is easy to see that solving \eqref{eqn:UQP} is equivalent to solving the linear system of equations
\begin{equation}
\label{eqn:equiLSE}
Px = q.
\end{equation}
Thus the unique optimal point to \eqref{eqn:UQP} is
\begin{equation}
\label{eqn:opt}
x_{\text{opt}} = P^{-1}q
\end{equation}
and the optimal value is
\begin{equation}
\label{eqn:optVal}
f(x_{\text{opt}}) = r - \frac{1}{2} q^t P^{-1} q.
\end{equation}
The matrix $P$ induces the norm
\begin{equation}
\label{eqn:NormP}
||x||_P := \sqrt{x^tPx}.
\end{equation}
In view of \eqref{eqn:UQP}, \eqref{eqn:optVal}, and \eqref{eqn:NormP}, we have
\begin{equation}
\label{eqn:fValNormPRel}
f(x) - f(x_{\text{opt}}) = \frac{1}{2}||x - x_{\text{opt}}||_P^2.
\end{equation}
A method that uses knowledge of $P,q$ and $r$ and obtains accurate estimates of \eqref{eqn:opt} and/or \eqref{eqn:optVal} is called a solution method for \eqref{eqn:UQP}.

We will say that the UQP in \eqref{eqn:UQP} is high dimensional, i.e., it shares the aforementioned high dimensional characteristics, if \emph{$n$ is high,} say of the order of $10^5$ or more, and \emph{the matrix $P$ is dense,} i.e., almost all its entries are nonzero. In what follows, we spell out the challenges involved in solving a high dimensional version of \eqref{eqn:UQP} when only limited computational resources are available.

\subsection{Computational Challenges for a High Dimensional UQP}
\label{subsec:Challenges}
\begin{itemize}
\item \textbf{There is a limit on the size of main memory and secondary storage:} Under present hardware technology, main memory (i.e., RAM) size of a processor is of the order of a few GigaBytes (GB) while secondary storage devices go up to few TeraBytes (TB). But observe, as an example, that when $n = 2^{20} \approx 10^6,$ the minimum space required just to store the $n^2$ entries of $P$ is $8$ TeraBytes---assuming double precision, i.e., $8$ Bytes per entry. This implies that to tackle a high dimensional UQP, we firstly need several secondary storage devices to store $P.$ Second, it becomes a necessity to slice the large $P$ into numerous chunks and store them in a distributed manner across multiple secondary storage devices. Third, only a few chunks, i.e., a small portion of the large $P,$ can be loaded into the main memory of a single processor at any given time.

\item \textbf{Computing $\mathbf{f(x)}$ or $\mathbf{\nabla f(x)}$ has long running time:} With the above memory issues, reading a large $P$ matrix into the main memory of a single processor, if there is such a need, can be done only in a chunk by chunk manner. This implies that, with a single processor, even the simple matrix vector multiplication $Px$ will have to be done only in a piecemeal fashion in high dimensions. Although the time complexity of this operation still remains $O(n^2),$ multiple secondary storage accesses, especially from noncontiguous locations, will add a nontrivial overhead.

    But observe from \eqref{eqn:UQP} and \eqref{eqn:gradUQP} that computing $f(x)$ or $\nabla f(x)$ at any arbitrary point $x \in \mathbb{R}^n$ necessarily involves the operation $Px.$ Consequently, in high dimensions, {\em any method that uses function values and/or gradient values to iteratively improve the solution quality will have to suffer an extremely large running time per iteration on a single processor;} several hours per iteration is not unimaginable for $n$ just a million in size.

\item \textbf{Parallelization offers little help:}  The running time of the operation $Px$ discussed above can definitely be cut down by parallelizing the computation across multiple processors. However, due to budgetary constraints, practitioners usually have access to only a fixed finite number of processors at any given time. These limited resources imply that, in high dimensions, the per iteration running time of methods that use function/gradient values to iteratively improve the solution quality can be cut down only marginally, usually only by a constant factor.
\end{itemize}

We would like to mention here that if \eqref{eqn:UQP} can be decomposed into a small number of independent problems then each of these problems can be accommodated and solved in parallel on independent machines. This idea is applicable whenever {\em the matrix $P$ has an inherent block diagonal structure with small number of blocks}. But our interest lies in high dimensional setups where the matrix $P$ is dense. As pointed out by \cite{richtarik2012parallel}, the only way here may be to appeal to a serial method leveraging a single processor or a fixed number of processors.

Circumventing the above challenges is the aim here. That is, we presume the availability of only a simplistic computational setup throughout the paper. This setup, referred to henceforth as the {\em finite main memory} (FMM) setup, is assumed to have:
\begin{enumerate}
\item a single processor with limited main memory and

\item finite sized secondary storage devices in numbers large enough to store the input data of the high dimensional UQP that is to be solved.
\end{enumerate}

Our goal then is to develop methods that can solve \eqref{eqn:UQP} even when $n$ is so large that the size of $P$ far exceeds the limited main memory size in the available FMM setup. Note that in theoretical analysis of a serial method having a single processor or finite number of processors makes no difference. For this reason and for pedagogical convenience, the FMM setup is assumed to have only a single processor. While discussing parallel and distributed implementation of the proposed GBCD method in Section~\ref{sec:ParDisImp}, we will use setups with multiple processors.

\subsection{Desiderata for a High Dimensional UQP Method}
\label{subsec:ProbDef}
We describe here a set of features that practitioners would prefer in an ideal solution method for a high dimensional version of \eqref{eqn:UQP} given only the FMM setup.

{\renewcommand*\theenumi{$\pmb{\mathcal{F}}_\arabic{enumi}$}
\begin{enumerate}
\item {\bf Work with only subset of entries of $P$ at one time.} The reason being that when $n$ is very large, the available main memory in FMM setup will be orders of magnitude smaller than the total data size of the matrix $P.$

\item {\bf Possess low per iteration running time (subquadratic in $n$).} In view of the discussion in the previous subsection, this feature essentially says that the method should never compute either $f(x)$ or $\nabla f(x)$ explicitly for any $x$ during the entire course of its execution.

\item {\bf Use hard partitioning of $P$.} That is, the method should suggest an explicit partition of $P$ into chunks that can be stored across multiple hard disks. After the partitioning, no major shuffling of the data should be required. This is needed because moving data within/across secondary hard disks is an extremely time consuming process. Furthermore, at any given time, the method should require to load data only from one or a few select chunks of $P$ into the main memory . The reason being that disk access time for non-contiguous data, as against contiguous data stored in the form of a chunk, is extremely high.
\end{enumerate}}

From now on, we will say that {\bf a solution method for \eqref{eqn:UQP} is high dimension compliant (HDC)} if and only if it has all the three features mentioned above. In addition to these basic features, one would optionally prefer that the HDC method has following desirable features as well.
{\renewcommand*\theenumi{$\pmb{\mathcal{F}}_\arabic{enumi}$}
\begin{enumerate}
\setcounter{enumi}{3}
\item {\bf Possess comparable running time.} That is, the total running time (iteration run time $\times$ number of iterations) of the method to find an approximate solution to $x_{\text{opt}}$ should be no worse than the best of the existing methods.

\item {\bf Give scope for parallel and distributed implementation.} That is, given a parallel and/or distributed computing environment, the method should be able to take advantage and achieve a speedup in its execution time by some factor.
\end{enumerate}}

\subsection{Contributions and Outline}
\label{subsec:OrgPaper}
In Section~\ref{subsec:priorArt}, we survey popular UQP solution methods and show that BK and BCD are the only existing methods that can be made HDC. In section~\ref{sec:UniView}, we prove that even natural generalizations of BK and BCD are not HDC. To the best of our knowledge, we are the first ones to prove such results. These results will show that coming up with better BK and BCD methods is indeed the right way forward for developing solution methods for high dimensional UQPs. In line with this view, we propose the GBCD method in Section~\ref{sec:GBCD}. We also discuss results here which show that the equivalent greedy BK method is non-HDC. We establish bounds on the convergence rate of GBCD method in Section~\ref{sec:PerMeas} and discuss heuristic ways to improve it in Section~\ref{sec:ParStra}. Parallel and distributed implementations of GBCD are discussed in Section~\ref{sec:ParDisImp}. Simulation results comparing GBCD with popular UQP solution methods are given in Section~\ref{sec:SimRes}. We finally conclude in Section~\ref{sec:Concl}.

\section{Survey of Existing UQP Methods and Research Gaps}
\label{subsec:priorArt}
We survey here the important UQP solution methods by grouping them into three categories---direct methods, classical iterative methods, and HDC methods. To solve a UQP, broadly speaking, the members of the first two categories essentially require that the input matrix be entirely available in the main memory. Consequently, these methods, unless modified, are inappropriate for solving high dimensional UQPs in an FMM setup. The third category includes the BK and BCD class of methods. As will be seen, these are the only existing methods that can be easily made HDC. After this survey, we discuss research gaps in state-of-the-art literature on BK and BCD methods.

\subsection{Direct Methods}
To solve the UQP in \eqref{eqn:UQP}, direct methods, see \cite{bertsekas1989parallel, davis2006direct}, begin by factorizing $P$ using techniques such as SVD, Cholesky decomposition, etc. The factored form of $P$ is then used to find \eqref{eqn:opt} and \eqref{eqn:optVal} in finite steps. These methods take $O(n^3)$ time, especially if $P$ is dense, before giving out the final solution. Popular members in this category include Gaussian elimination, Cholesky factorization, etc. When $n$ is really small, say less than $10^3,$ direct methods are highly preferred. However, the lack of intermediate solutions coupled with the enormous difficulty in factorizing a matrix, when it is stored in a distributed fashion, makes direct methods unfavourable for solving high dimensional UQPs in the FMM setup. These methods are clearly non-HDC.

\subsection{Classical Iterative Methods}
\label{subsec:classIterMethods}
These methods solve a UQP by using its function/gradient values repeatedly to generate a sequence of improving solution estimates of \eqref{eqn:opt} and/or \eqref{eqn:optVal}. These are strongly favored when $n$ is moderately large, say roughly from $10^3$ upto $10^5.$ However, the $O(n^2)$ running time per iteration for these methods makes them non-HDC. Specifically, these methods tend to have exorbitant waiting times before each update of the solution estimate when used to solve a high dimensional UQP. If, however, this is not an issue, then, after appropriate modifications, one can go ahead and use these methods to solve high dimensional UQPs even in the FMM setup.

Classical iterative methods for UQPs can be broadly subdivided into line search methods, cutting plane methods, and direct search methods. A line search method, see \cite{wright1999numerical}, uses the following idea in every iteration. It first decides upon a search direction. Along this direction starting from the current estimate, the point that is closest to the optimal point of the given UQP, i.e., $x_{\text{opt}},$ under appropriate norm, is then declared the new estimate. Important algorithms here include the steepest descent and the conjugate gradient. Both these methods use the gradient, after suitable modifications, to generate the search direction in each iteration. The $||\cdot||_P$ norm from \eqref{eqn:NormP} is used to decide the new estimate in each iteration.

Another category of UQP methods that use gradients are the cutting plane methods, see \cite{kelley1960cutting, levin1965algorithm}. Here, the idea is to start with a feasible set containing $x_{\text{opt}}.$ Using a hyperplane normal to the gradient at its centre of gravity, the feasible set is cut into two parts. The side containing $x_{\text{opt}}$ is then called the new feasible set and the above idea is repeated. A major difficulty here is the efficient computation of the centre of gravity of the feasible set. This can itself turn out into a very hard problem if the feasible region is arbitrary. Alternatives such as the ellipsoid method, see \cite{bland1981ellipsoid}, etc., have been proposed to overcome this issue. However, these alternatives need suitable revisions before they can be used to solve a high dimensional UQP in the FMM setup.

Direct search methods, see \cite{kolda2003optimization}, avoid gradients and use only function values to sequentially improve the solution quality. The \cite{nelder1965simplex} simplex algorithm is a key member here. To solve a UQP, this method begins with a $n-$dimensional simplex. Using a sequence of  reflection, expansion, contraction, and reduction type transformations on the simplex, the method sequentially reduces the function values at each of its vertices. The Nelder-Mead approach is, however, known to be slower in convergence in comparison to gradient based schemes described above. Hence this simplex based approach should not be used to solve a high dimensional UQP in the FMM setup even if $O(n^2)$ running time per iteration is not an issue.

\subsection{HDC Methods}
\label{subsec:HDCMethods}
We include in this category exclusively the class of BK and BCD methods. These are also iterative methods that follow the line search principle. However, relative to classical line search methods mentioned before, the BK and BCD methods use a very different approach to solve a UQP. First, the BK and BCD methods never explicitly compute function/gradient values. Second, in each iteration, they need only finite rows of the input matrix. As we shall see at the end of this subsection, these are also the reasons why BK and BCD methods can be made HDC.

\subsubsection{Block Kaczmarz (BK) Methods:}
The fundamental member of this class is the simple \cite{kaczmarz1937angenaherte} method. To solve the UQP in \eqref{eqn:UQP}, this method uses precisely one row of the input matrix $P$ in every iteration. Its update rule is given by
\begin{equation}
\label{eqn:updRuleSK}
x_{k + 1} = x_k + P_{i}^t\;(P_{i}P_{i}^t)^{-1} \;(q_i - P_{i}x_k),
\end{equation}
where $P_{i}$ is the $i^{th}$ row of matrix $P.$ The index $i$ is chosen in a round robin fashion from one iteration to another.

\cite{elfving1980block} generalized the above idea. We will refer to his method as the traditional BK. As opposed to \eqref{eqn:updRuleSK}, the traditional BK method uses multiple rows of $P$ per iteration. More formally, the traditional BK method first partitions the rows of $P$ into different blocks. Suppose that the row submatrices $P_{\pi_1}, \ldots, P_{\pi_m}$ denote these blocks, where $\Pi \equiv \{\pi_1, \ldots, \pi_m\}$ represents a partition of the row index set $[n]:= \{1, \ldots, n\}.$ The traditional BK method then uses these blocks in a round robin fashion across iterations for improving the solution estimate. In notations, its update rule is
\begin{equation}
\label{eqn:updRuleBK}
x_{k+1} = x_k + P_{\pi_i}^t \; (P_{\pi_i}P_{\pi_i}^t)^{-1} \; (q_{\pi_i} - P_{\pi_i}x_k),
\end{equation}
where $q_{\pi_i}$ is the subvector of $q$ corresponding to the row indices given in $\pi_i$.

Clearly, there are two ways in which the traditional BK method can be modified.
\begin{enumerate}
\item {\bf Input partition $\Pi$:} Which partition of the rows of $P$ to use at the beginning of the method?

\item {\bf Block Selection Strategy $\zeta_b$:} Which block to work with in each iteration of the method?
\end{enumerate}

By choosing appropriately these strategies, one can obtain different BK methods. The BK method in generic form is given in Algorithm~\ref{alg:BKSch}. We emphasize here that the convergence rate of a BK method strongly depends on the choice for both $\Pi$ and $\zeta_b.$

\begin{algorithm}[ht!]
   \caption{Generic Block Kaczmarz (BK) Method}
   \label{alg:BKSch}
\begin{algorithmic}
    \STATE {\bfseries Input:} initial guess $x_0 \in \mathbb{R}^n,$ partition $\Pi \equiv \{\pi_1, \ldots, \pi_m\}$ and strategy $\zeta_b$
    \STATE {\bfseries Main Procedure:}
    \FOR{$k = 0, 1, \ldots,$}
    \STATE Using strategy $\zeta_b$, pick a block $P_{\pi_i}$ from the available set of blocks.
    \STATE $x_{k + 1} \leftarrow x_k + P_{\pi_i}^t (P_{\pi_i} P_{\pi_i}^t)^{-1} (q_{\pi_i} - P_{\pi_i} x_k).$
    \ENDFOR
\end{algorithmic}
\end{algorithm}

\subsubsection{Block Coordinate Descent (BCD) Methods:} The counterpart of the simple Kaczmarz method in the BCD class is the Gauss-Seidel method, see \cite{golub2012matrix}. To solve the UQP in \eqref{eqn:UQP}, this also selects one row of $P$ per iteration using round robin. However, unlike the simple Kaczmarz, this method updates only one coordinate of the estimate---index of which matches the index of the chosen row of $P$---in every iteration. We will refer to the counterpart of the traditional BK method as the traditional BCD method, see \cite{hildreth1957quadratic} and \cite{warga1963minimizing}. In this method also, rows of $P$ are first partitioned into arbitrary blocks $P_{\pi_1}, \ldots, P_{\pi_m}.$ In each iteration, the traditional BCD method then chooses one block in a round robin fashion and improves the solution estimate at only those coordinates whose indices match with those of the chosen rows of $P.$ Once again, by choosing a partition $\Pi$ of $P$ and a strategy $\zeta_b$ to select the appropriate block in each iteration, one can come up with different BCD methods. The generic BCD method is given in Algorithm~\ref{alg:BCDSch}. Note in the update step that $\mathbb{I}_{\pi_i}$ denotes the rows of the $n \times n$ identity matrix with indices in $\pi_i$ and  $P_{\pi_i \pi_i} := P_{\pi_i} \mathbb{I}_{\pi_i}^t.$ Like BK methods, the convergence rate of the BCD methods strongly depends on the choice for both $\Pi$ and $\zeta_b$.

\begin{algorithm}[ht!]
   \caption{Generic Block Coordinate Descent (BCD) Method}
   \label{alg:BCDSch}
\begin{algorithmic}
    \STATE {\bfseries Input:} initial guess $x_0 \in \mathbb{R}^n,$ partition $\Pi \equiv \{\pi_1, \ldots, \pi_m\}$ and strategy $\zeta_b$
    \STATE {\bfseries Main Procedure:}
    \FOR{$k = 0, 1, \ldots,$}
    \STATE Using strategy $\zeta_b$, pick a block $P_{\pi_i}$ from the available set of blocks.
    \STATE $x_{k + 1} \leftarrow x_k + {\mathbb{I}}_{\pi_i}^t (P_{\pi_i \pi_i})^{-1} (q_{\pi_i} - P_{\pi_i} x_k).$
    \ENDFOR
\end{algorithmic}
\end{algorithm}

We now make some remarks concerning the BK and BCD methods. This will show how these methods can be made HDC. For a UQP, let $\rho(n)$ denote the maximum number of rows of $P$ that can be stored in the main memory of the given FMM setup. For a partition $\Pi$ of $P,$ let $d$ denote the maximum of number of rows in the blocks $P_{\pi_1}, \ldots, P_{\pi_m}.$

\begin{itemize}
\item For the update step in Algorithms~\ref{alg:BKSch} and \ref{alg:BCDSch}, at most $d$ rows of $P$ and the corresponding entries of $q$ are needed in the main memory. This translates to $O(nd)$ space.

\item The update step in BK method has a running time of $O(nd^2).$ This is subquadratic if $d < \sqrt{n}.$

\item Similarly, for the update step in BCD method, the running time is $O(nd) + O(d^3)$ which is subquadratic if $d < n^{2/3}.$

\item Thus, choosing a strategy $\zeta_B$ that respects above space and time complexity bounds in each iteration and picking $d$ that is less than $\min\{\rho(n), \sqrt{n}\}$ for BK method and less than $\min\{\rho(n), n^{2/3}\}$ for BCD method will ensure that Algorithms~\ref{alg:BKSch} and \ref{alg:BCDSch} have features $\mathcal{F}_1$ and $\mathcal{F}_2.$

\item In both BK and BCD method, $P$ is a priori partitioned into blocks (row submatrices). These blocks can be stored in separate locations in the secondary storage and only one of them is needed in each iteration. This is clearly as desired in feature $\mathcal{F}_3.$ This verifies that the class of BK and BCD methods are the only existing methods that can be made HDC. This also shows the convenience of using BK and BCD methods for solving high dimensional UQPs in the FMM setup.
\end{itemize}

\subsection{Research Gaps in BK and BCD Methods}
\label{subsec:RGap}
As mentioned earlier, the convergence rate of BK and BCD methods crucially depend on the input partition $\Pi$ as well as block selection strategy $\zeta_b.$ It is however difficult to jointly analyze the effect of both $\Pi$ and $\zeta_b$ on the convergence rate. Hence, most of the existing literature first fixes an arbitrary partition of $P$ and then identifies a good choice for strategy $\zeta_b.$ Once this is done, effort is made to come up with a partition $\Pi$ that works best with the earlier chosen $\zeta_b.$

In line with the above viewpoint, fix an arbitrary row partition of $P.$ In traditional BK and BCD methods, the strategy $\zeta_b$ is to sweep through the given blocks in a round-robin manner across iterations. For these methods, \cite{deutsch1985rate, deutsch1997rate, galantai2005rate} have shown that it is not easy to even evaluate the convergence rate in typical matrix quantities, let alone finding the associated best row partition of $P.$ To overcome this, recent works have suggested randomly picking blocks in each iteration. For example, randomized Kaczmarz algorithm of \cite{strohmer2009randomized} picks one row of $P$ in each iteration using a probability distribution that is proportional to the square of length of its rows. Randomized BK method from \cite{needell2014paved}, dealing with multiple rows of $P$ per iteration, picks blocks in a uniformly random fashion. Here, however, each row of $P$ needs to have unit $||\cdot||_2$ norm. The randomized coordinate descent method of \cite{leventhal2010randomized} works by picking a single row in each iteration using a distribution that is proportional to the diagonal entries of $P.$ Randomized BCD method from \cite{nesterov2012efficiency} uses the eigenvalues of block diagonal submatrices of $P$ to select blocks in every iteration. In all these approaches, neat convergence rate estimates have been obtained and effort is now underway to find the corresponding best partition, see \cite{needell2014paved}.

It is worth mentioning that in most of the existing randomized BK and BCD methods including the ones mentioned above, the probability distribution for selecting the blocks in each iteration remains fixed throughout the execution. This impedes performance in several instances. For an example, see Experiment 2 of Section~\ref{sec:SimRes}. Adaptive block selection strategies are thus better. Adaptive versions of randomized coordinate descent algorithms have been discussed in \cite{loshchilov2011adaptive, glasmachers2013accelerated}. In this paper, our focus is on greedy block selection strategies which are adaptive but deterministic. Observe from Algorithms~\ref{alg:BKSch} and \ref{alg:BCDSch} that, in each iteration, the new estimate can be generated using any of the available blocks of $P.$ Hence a natural greedy block selection strategy is to pick the best block in each iteration. That is, pick the block amongst all possible candidate choices for which the revised estimate $x_{k+1}$ is the closest under an appropriate norm to $x_{\text{opt}}.$ The BK method with this greedy strategy, however, ends up having a per iteration run time of $O(n^2)$ and is consequently non-HDC. This fact is briefly discussed at the end of Section~\ref{sec:GBCD}.  An alternative to this greedy BK approach via lower dimensional projections has been given by \cite{eldar2011acceleration}. In this paper, we study the BCD method with the above greedy approach.

All BK and BCD methods described earlier including the proposed GBCD method take in as input a hard partition of the $P$ matrix and work with only one of the blocks in each iteration. Note that it is possible to alleviate this restriction and develop methods that are allowed to pick any arbitrary rows of $P$ in an iteration. The BCD method proposed by \cite{bo_uai08} is in fact of this kind. We will refer to this method as GBCD-BS and it works as follows. It first fixes a number $r < n.$ Then in each iteration, it greedily picks $r$ rows which will help descend the largest in that iteration. We wish to emphasize here that the GBCD-BS method tries to identify the best $r$ rows in each iteration while our proposed GBCD method tries to identify the best block of rows from a prefixed partition. Ideally speaking, when both methods work with roughly the same number of rows in each iteration, the GBCD-BS method should outperform the proposed GBCD method. However, when solving high dimensional UQPs using only the FMM setup, the proposed GBCD method will typically be better. The reason being the secondary storage fetches that will be needed in each and every iteration. For the proposed GBCD method, the rows will have to be fetched from contiguous locations while for the GBCD-BS method the rows will typically be from non-contiguous locations. Hence, the per iteration running time of the proposed GBCD method will be significantly lower than that of the GBCD-BS giving it the above mentioned advantage. A demonstration of this has been given in Experiment 1 of Section~\ref{sec:SimRes}. We wish to mention here that when both the GBCD-BS as well as the proposed GBCD method work with exactly one row in each iteration, then they are one and the same.

A comparative analysis of the proposed GBCD method against recent popular approaches for UQP is given in Table~\ref{tab:compAnalysis}. We have skipped the GBCD-BS method as its theoretical convergence rate estimate is not available. Note that the error definition is different for different methods. Although it is possible to express the estimates of the table using a common error definition, we refrain from doing so. This is because conversion inequalities will add approximation errors leading to an unfair comparison. A summary of the notations used in the table are as follows. By an $\epsilon-$close solution for a method, we imply an estimate $x_k$ for which the corresponding $k-$step error is less than $\epsilon^2.$ The time for $\epsilon-$close solution is the product of iteration run time and number of iterations required to find an $\epsilon-$close solution. We have assumed that the BK and BCD methods---except randomized Kaczmarz and randomized coordinate descent---use a partition of $P$ that has $m$ blocks and each block is made up at most $d$ rows. The maximum and minimum eigenvalue functions are denoted using $\lambda_{\max}(\cdot)$ and $\lambda_{min}(\cdot).$ The usual condition number of $P$ is denoted by $\kappa(P):= \lambda_{\max}(P)/\lambda_{\min}(P).$ Related to this is the scaled condition number, introduced by \cite{demmel1988probability}, which is defined as $\tilde{\kappa}(P) := ||P||_F/\lambda_{\min}(P),$ where $||\cdot||_F$ denotes the Frobenius norm. The trace of $P$ is denoted using $\text{Tr}(P).$ We use $\Pi$ to denote a partition of the row index set $[n]$ of $P.$ Further, $\beta := \max_{\pi \in \Pi} \lambda_{\max} (P_{\pi \cdot} P_{\pi \cdot}^t).$ Finally, $P_{\Pi}$ and $B_{\Pi}$ are as defined in \eqref{eqn:PermP} and \eqref{eqn:BlkDiag}.

\begin{sidewaystable}
\centering
\vspace{41.5em}
\caption{\label{tab:compAnalysis} Comparison of Proposed GBCD method with Popular UQP methods.}
\begin{tabular}{ l | c | c | c | c }
\hline
 Method & \begin{tabular}{@{}c@{}}Iteration \\ Run Time \end{tabular} & \begin{tabular}{@{}c@{}}$k-$step Error \\ Definition \end{tabular} & \begin{tabular}{@{}c@{}} Upper Bound on \\ $k-$step Error\end{tabular} & \begin{tabular}{@{}c@{}} Upper Bound on Time \\ for $\epsilon-$close solution \end{tabular} \\
\hline
 & & & & \\
Steepest Descent & $O(n^2)$ & $\dfrac{||x_k - x_{\text{opt}}||^2_P}{||x_0 - x_{\text{opt}}||_P^2}$ & $\left[\dfrac{\kappa(P) - 1}{\kappa(P) + 1}\right]^{2k}$ & $O(n^2 \kappa(P) \log(\frac{1}{\epsilon}))$\\
 & & & & \\
Conjugate Gradient & $O(n^2)$ & $\dfrac{||x_k - x_{\text{opt}}||^2_P}{||x_0 - x_{\text{opt}}||_P^2}$ & $4\left[\dfrac{\sqrt{\kappa(P)} - 1}{\sqrt{\kappa(P)} + 1}\right]^{2k}$ & $O(n^2 \sqrt{\kappa(P)} \log(\frac{1}{\epsilon}))$ \\
 & & & & \\
\begin{tabular}{@{}l@{}}Randomized Kaczmarz \\ -\cite{strohmer2009randomized}\end{tabular} & $O(n)$ & $\mathbb{E}\left[\dfrac{||x_k - x_{\text{opt}}||^2_2}{||x_0 - x_{\text{opt}}||_2^2}\right]$ & $\left[1 - \dfrac{1}{\tilde{\kappa}(P)^2}\right]^k$ & $O(n \tilde{\kappa}(P)^2 \log(\frac{1}{\epsilon}))$ \\
 & & & & \\
\begin{tabular}{@{}l@{}}Randomized BK \\ -\cite{needell2014paved} \end{tabular} & $O(nd^2)$ & $\mathbb{E}\left[\dfrac{||x_k - x_{\text{opt}}||^2_2}{||x_0 - x_{\text{opt}}||_2^2}\right]$ & $\left[1 - \dfrac{\lambda_{\min}^2(P)}{m \beta}\right]^k$ & $O\left(\dfrac{nm \beta \log(\frac{1}{\epsilon})}{\lambda_{\min}^2(P)}\right)$\\
 & & & & \\
 \begin{tabular}{@{}l@{}}Randomized Coordinate Descent \\ -\cite{leventhal2010randomized} \end{tabular}& $O(n)$ & $\mathbb{E}\left[\dfrac{||x_k - x_{\text{opt}}||^2_P}{||x_0 - x_{\text{opt}}||_P^2}\right]$ & $\left[1 - \dfrac{\lambda_{\min}(P)}{\text{Tr}(P)}\right]^k $ & $O\left(\dfrac{n \text{Tr}(P)\log(\frac{1}{\epsilon})}{ \lambda_{\min}(P)}\right)$\\
 & & & & \\
\begin{tabular}{@{}l@{}}Randomized BCD  ($d < \sqrt{n}$) \\ -\cite{nesterov2012efficiency} \end{tabular} & $O(nd)$ & $\mathbb{E}\left[\dfrac{||x_k - x_{\text{opt}}||^2_P}{||x_0 - x_{\text{opt}}||_P^2}\right]$ & $\left[1 - \dfrac{\lambda_{\min}(P)}{\sum\limits_{\pi \in \Pi} \lambda_{\max}(P_{\pi\pi})}\right]^k$ & $O\left(\dfrac{n \sum\limits_{\pi \in \Pi} \lambda_{\max}(P_{\pi\pi}) \log(\frac{1}{\epsilon})}{ \lambda_{\min}(P)}\right)$\\
 & & & & \\
\begin{tabular}{@{}l@{}}GBCD  ($d < \sqrt{n}$) \\ -Proposed Method \end{tabular} & $O(nd)$ & $\dfrac{||x_k - x_{\text{opt}}||^2_P}{||x_0 - x_{\text{opt}}||_P^2}$ & $\left[1 - \dfrac{\lambda_{\min}(P_{\Pi} B_{\Pi}^{-1})}{m} \right]^k$ & $O\left(\dfrac{nmd \log(\frac{1}{\epsilon})}{\lambda_{\min}(P_{\Pi} B_{\Pi}^{-1})}\right)$\\
 & & & & \\
\hline
\end{tabular}
\end{sidewaystable}

\section{Are There More HDC Methods?}
\label{sec:UniView}
In the previous section, we showed that the class of BK and BCD are the only existing HDC methods for UQP. This hints at the possibility that by generalizing the BK and BCD philosophies, one may be able to come up with a richer class of HDC solution methods for UQP. To check this out, we first propose the framework of descent via lower dimensional restrictions (DLDR)---a generic approach to solve high dimensional convex optimization programs (not necessarily quadratic)---in Subsection~\ref{subsec:DLDR}. In Subsections~\ref{subsec:BK} and \ref{subsec:BCD}, we then show respectively that BK and BCD methods for UQP are specialized instances of DLDR. Using this, we finally conclude that  even natural generalizations of BK and BCD are not HDC.

\subsection{Descent via Lower Dimensional Restrictions (DLDR) Framework}
\label{subsec:DLDR}
Consider the unconstrained convex optimization program
\begin{equation}
\label{eqn:optProgram}
\underset{x \in \mathbb{R}^n}{\text{min}} \hspace{0.5em} g(x),
\end{equation}
where $n$ is large and $g:\mathbb{R}^n \rightarrow \mathbb{R}$ is a strict convex function bounded from below. To solve \eqref{eqn:optProgram}, we suggest the following iterative idea. \emph{Given the current estimate, pick an appropriate affine space passing through it. Declare the optimum of $g$ restricted to this affine space as the new estimate.} Note that, as $g$ is a strict convex function, its affine restriction will again be strictly convex but of lower dimensions. The above idea forms the basis of proposed DLDR framework. Details are given in Algorithm~\ref{alg:DLDR}. The notation $\text{col}(M_k)$ denotes the column space of matrix $M_k$ and $\ell.i.$ stands for linearly independent. When $d_k = 1$ for each $k,$ the DLDR framework is precisely the classical line search method for optimization. The DLDR framework is thus its canonical generalization.

\begin{algorithm}[ht!]
   \caption{DLDR Framework}
   \label{alg:DLDR}
\begin{algorithmic}
    \STATE {\bfseries Input:} initial guess $x_0 \in \mathbb{R}^n,$ strategy $\zeta.$
    \STATE {\bfseries Main Procedure:}
    \FOR{$k = 0, 1, \ldots,$}
    \STATE Choose $M_k \in \mathbb{R}^{n \times d_k},$ $d_k < n,$ with $\ell.i.$ columns using strategy $\zeta.$
    \STATE Define affine space $\mathcal{A}_k$ as $x_k + \text{col}(M_k).$
    \STATE $x_{k + 1} \leftarrow \underset{x \in \mathcal{A}_k}{\text{argmin}} \hspace{0.5em} g(x).$
    \ENDFOR
\end{algorithmic}
\end{algorithm}

The following fact concerning the DLDR framework is easy to see.

\begin{fact} [Guaranteed descent in each iteration] \label{fact:GuaDesc}
Whatever be the strategy $\zeta$ in DLDR framework, for each $k$, we always have
\begin{equation}
\label{eqn:gtdDes}
g(x_{k + 1}) \leq g(x_k)
\end{equation}
\end{fact}

\subsection{BK Method and DLDR Framework}
\label{subsec:BK}
We first describe how the DLDR framework can be used to solve the UQP in \eqref{eqn:UQP}. Using $x_{\text{opt}}$ given in \eqref{eqn:opt}, define for each $x \in \mathbb{R}^n,$
\begin{equation}
\label{eqn:g_BK}
g(x) = ||x - x_{\text{opt}}||_2^2.
\end{equation}
Clearly, $g$ is a strict convex function having $x_{\text{opt}}$ as its unique optimal point. Hence, it follows that one can solve the UQP in \eqref{eqn:UQP} by alternatively finding the minima of \eqref{eqn:g_BK} using the DLDR framework. The following fact is now immediate.

\begin{fact}
\label{fact:twonorm}
The update rule of the DLDR framework for \eqref{eqn:g_BK} is:
\begin{equation}
\label{eqn:BKupd}
x_{k + 1} = x_k + M_k (M_k^t M_k)^{-1} M_k^t (x_{\text{opt}} - x_k).
\end{equation}
\end{fact}

This fact establishes the desired connection between BK method and DLDR framework .

\begin{lemma}
\label{lem:BKisDLDR}
Consider the BK method from Algorithm~\ref{alg:BKSch} and the DLDR framework from Algorithm~\ref{alg:DLDR}. Suppose that $\zeta = (\Pi, \zeta_b),$ i.e., if in the $k-$th iteration $\zeta_b$ suggests picking the block $P_{\pi_i},$ then $\zeta$ sets $M_k = P_{\pi_i}^t.$ Then, applying the BK method to solve the UQP in \eqref{eqn:UQP} is precisely the same as if we were applying the DLDR framework to find the minima of \eqref{eqn:g_BK}.
\end{lemma}

From Fact~\ref{fact:GuaDesc}, we know that there is guaranteed descent in each iteration of the DLDR framework even if the matrix sequence $\{M_k\}$ is chosen completely arbitrarily. This may appear to suggest that it is possible to generalize the idea of the BK method and obtain additional HDC solution methods for UQP. The following result shows that this is not true.

\begin{lemma}
\label{lem:ImpDLDRKacz}
If DLDR framework is used to find the minima of \eqref{eqn:g_BK}, then it is possible to implement it in practice only when each of the matrix $M_k$ is generated using the rows of $P$.
\end{lemma}

This follows from the fact that the vector $M_k^t x_{\text{opt}}$ in the update rule of DLDR, see \eqref{eqn:BKupd},  would be unknown if the matrix $M_k$ is arbitrary and it is only for the special case of $M_k = P_{\pi_i}^t,$ $\pi_i \in \Pi,$ where $M_k^t x_{\text{opt}}$ can be replaced with $q_{\pi_i}.$ This result essentially says that BK method from Algorithm~\ref{alg:BKSch} is the only possible form of the DLDR framework when finding the minima of \eqref{eqn:g_BK}. That is, no natural generalization of the BK method exists; let alone HDC versions.

\subsection{BCD Method and DLDR Framework}
\label{subsec:BCD}

Let $x_{\text{opt}}$ be as in \eqref{eqn:opt}. In contrast to \eqref{eqn:g_BK}, let
\begin{equation}
\label{eqn:g_BCD}
g(x) = ||x - x_{\text{opt}}||_P^2
\end{equation}
Clearly, $g$ is again strictly convex and $x_{\text{opt}}$ is its unique optimal point. Once again, minimizing \eqref{eqn:g_BCD} using the DLDR framework will help solve the UQP in \eqref{eqn:UQP}. The following fact is easy to see.

\begin{fact}
\label{fact:Pnorm}
The update rule of the DLDR framework for \eqref{eqn:g_BCD} is:
\begin{equation}
\label{eqn:BCDupd}
x_{k + 1} = x_k + M_k (M_k^tPM_k)^{-1} M_k^t P(x_{\text{opt}} - x_k).
\end{equation}
\end{fact}

This gives us the following desired result.
\begin{lemma}
\label{lem:BCDisDLDR}
Consider the BCD method from Algorithm~\ref{alg:BCDSch} and the DLDR framework from Algorithm~\ref{alg:DLDR}. Suppose that $\zeta = (\Pi, \zeta_b),$ i.e., if in the $k-$th iteration $\zeta_b$ suggests picking the block $P_{\pi_i},$ then $\zeta$ sets $M_k = \mathbb{I}_{\pi_i}^t.$ Then, applying the BCD method to solve the UQP in \eqref{eqn:UQP} is precisely the same as if we were applying the DLDR framework to find the minima of \eqref{eqn:g_BCD}.
\end{lemma}

Using \eqref{eqn:opt}, observe that \eqref{eqn:BCDupd} can be rewritten as
\begin{equation}
\label{eqn:BCDupdMod}
x_{k + 1} = x_k + M_k (M_k^tPM_k)^{-1} (M_k^t q - M_k^tPx_k).
\end{equation}
Because of this, we have the following result which is in complete contrast to Lemma~\ref{lem:ImpDLDRKacz}.

\begin{lemma}
\label{eqn:ImpDLDRBCD}
The DLDR framework of Algorithm~\ref{alg:DLDR} for \eqref{eqn:g_BCD} is implementable whatever be the choice for the sequence of full column rank matrices $\{M_k\}_{k \geq 0}.$
\end{lemma}

Thus, the BCD idea can indeed be generalized in numerous ways. But observe from \eqref{eqn:BCDupdMod} that even if one column of $M_k$ is dense then the time required to compute $M_k^tP,$ and hence the per iteration running time, will be $O(n^2).$ In fact, majority of the rows of $M_k$ must be all zero vectors. Otherwise we will need to store a lot of entries of $P$ in the main memory. Hence it follows that BCD from Algorithm~\ref{alg:BCDSch}, other than minor modifications, is the only HDC form of the DLDR framework when finding the minima of \eqref{eqn:g_BCD}.

The discussion in this and the previous section confirms the following: \emph{Coming up with better BK and BCD methods is indeed the right way forward in developing solution methods for high dimensional UQPs.} In line with this view, we propose the GBCD method in the next section. The greedy BK method unfortunately has a per iteration run time of $O(n^2)$ and hence is not HDC. The reasons for this are also briefly discussed at the end of the next section.

\section{Greedy Block Coordinate Descent (GBCD) Method}
\label{sec:GBCD}
Consider the generic BCD method given in Algorithm~\ref{alg:BCDSch}. As we had mentioned in Subsection~\ref{subsec:RGap}, the idea of greedy block selection strategy is to pick, in some sense, the best block in each iteration. Our first goal here is to understand this notion of best block in the context of BCD methods. From Lemma~\ref{lem:BCDisDLDR} and Fact \ref{fact:GuaDesc} the following result is immediate.
\begin{lemma}
\label{lem:BCDMonCov}
Let $\{x_k\}$ be the estimates generated by the BCD method of Algorithm~\ref{alg:BCDSch}. Then, whatever be the partition $\Pi$ and strategy $\zeta_b$, $||x_k - x_{\text{opt}}||_P$ is a non-increasing function of $k.$
\end{lemma}

In fact, one can find the exact relation between $||x_k - x_{\text{opt}}||_P$ and $||x_{k + 1} - x_{\text{opt}}||_P.$ For $\pi_i \in \Pi,$ let
\begin{equation}
\label{eqn:partGrad}
\nabla_{\pi_i}f(x) := P_{\pi_i} x - q_{\pi_i}
\end{equation}
denote the partial gradient vector and let
\begin{equation}
\label{eqn:beta}
\beta_{\pi_i}(x) := \nabla_{\pi_i} f(x)^t P_{\pi_i\pi_i}^{-1} \nabla_{\pi_i} f(x),
 \end{equation}
where $P_{\pi_i}$ and $P_{\pi_i\pi_i}$ are as in Algorithm~\ref{alg:BCDSch}. Note that $\beta_{\pi_i}(\cdot)$ is non-negative function for each $\pi_i.$

\begin{lemma}
\label{lem:PythThmNrmP}
Suppose that strategy $\zeta_b$ in Algorithm~\ref{alg:BCDSch} suggests choosing block $P_{\pi_i}$ in iteration $k.$ Then,
\begin{equation}
\label{eqn:PythThmNrmP}
||x_{k + 1} - x_{\text{opt}}||_P^2 = ||x_k - x_{\text{opt}}||_P^2 - \beta_{\pi_i}(x_k).
\end{equation}
\end{lemma}
\proof{Proof.}
This is immediate.
\endproof

Based on this, we come up with the following definition for the best block.
\begin{definition}
In iteration $k$ of a BCD method, we will say $P_{\pi_j}$ is the best block amongst $P_{\pi_1}, \ldots, P_{\pi_m}$ if $\beta_{\pi_j}(x_k) \geq \beta_{\pi_i}(x_k)$ for each $i \neq j.$
\end{definition}

Note that choosing the best block will ensure that the revised estimate $x_{k + 1}$ is closest possible to $x_{\text{opt}}$ in $||\cdot||_P$ norm amongst available choices in iteration $k.$

The idea of the GBCD method that we propose is to select such a best block in every iteration. But observe that this strategy would require computing $\beta_{\pi_i}(x_k)$ for each $\pi_i$ in the $k-$th iteration. From \eqref{eqn:beta}, this in fact implies that one would have to compute the complete gradient $\nabla f(x_k)$ in the $k-$th iteration. This is a cause for serious concern as explicitly computing the gradient in each iteration would essentially make the GBCD method non-HDC. Fortunately, as the next result shows, computing gradients in BCD methods is very easy. This result follows mainly due to the fact that successive iterates in BCD methods, unlike other UQP solution methods, differ only in few coordinates. Hence, unlike for other methods, the running time to compute a gradient in BCD methods is significantly lower than $O(n^2).$

\begin{lemma}
\label{lem:gradUpd}
Suppose that strategy $\zeta_b$ in Algorithm~\ref{alg:BCDSch} suggests choosing block $P_{\pi_i}$ in iteration $k.$ Then,
\[
\nabla f(x_{k + 1}) = \nabla f(x_k) + P^t_{\pi_i} (P_{\pi_i \pi_i})^{-1} (q_{\pi_i} - P_{\pi_i} x_k).
\]
\end{lemma}
\proof{Proof.}
This follows from \eqref{eqn:gradUQP}, the update rule in Algorithm~\ref{alg:BCDSch} and the fact that $P \in S_{++}^n.$ \Halmos
\endproof

This result tells us that the gradient in each iteration can be computed in an iterative fashion. The advantage of doing so is that, once we know $\nabla f(x_0),$ computing the gradient in each iteration requires only $O(nd)$ time, where recall $d$ is the maximum of number of rows in the blocks $P_{\pi_1}, \ldots, P_{\pi_m}.$ This is good because if the inverses $P_{\pi_1 \pi_1}^{-1}, \ldots, P_{\pi_m \pi_m}^{-1}$ are precomputed, then finding all of $\beta_{\pi_1}(x_k), \ldots, \beta_{\pi_m \pi_m}(x_k)$ in each iteration $k$ is now only an $O(nd)$ operation. The only bottleneck that remains is computing $\nabla f(x_0)$ and the inverses $P_{\pi_1 \pi_1}^{-1}, \ldots, P_{\pi_m \pi_m}^{-1}.$ Clearly, if $x_0 = 0,$ then $\nabla f(x_0) = -q.$ That is, computing the initial gradient in the special case of the initial estimate being the origin requires only $O(n)$ time. In all other cases, computing the initial gradient is an $O(n^2)$ operation. Computing the $m$ inverses requires $O(nd^2)$ time. This is subquadratic if $d < \sqrt{n}.$ The $O(n^2)$ and $O(nd^2)$ running time is large no doubt. But these operations are only one time and their running times will get amortized over iterates. Hence, the per iteration time complexity computation of the GBCD method need not consider the time to find the initial gradient and the $m$ inverses. Also, note that computing $\nabla f(x_{k+1})$ from $\nabla f(x_k)$ requires the same rows of $P$ and the corresponding entries of $q$ that are required in the update step. That is, the main memory space requirement remains at $O(nd).$

The above discussions put together shows that the GBCD method can indeed be made HDC. The only requirement is that for the input partition $\Pi,$ $d < \min\{\sqrt{n}, \rho(n)\},$ where $\rho(n)$ is as defined at the end of Subsection~\ref{subsec:HDCMethods}. We will denote the set of such partitions as $\wp.$ Algorithm~\ref{alg:BCDSchGreedy} describes the proposed GBCD method to solve a high dimensional UQP given only the FMM setup. Note that we have explicitly mentioned a step that deals with loading $P_{\pi}$ and $q_{\pi}.$ This is done to highlight the fact that all entries of $P$ and $q$ needed in an iteration will have to be fetched from the secondary storage only during that iteration.

\begin{algorithm}[ht!]
   \caption{Greedy Block Coordinate Descent (GBCD) method}
   \label{alg:BCDSchGreedy}
\begin{algorithmic}
    \STATE {\bfseries Input:} initial guess $x_0 \in \mathbb{R}^n,$ partition $\Pi \equiv \{\pi_1, \ldots, \pi_{m}\} \in \wp$
    \STATE {\bfseries Preprocessing:}
    \STATE Store $P_{\pi_i}$ and $q_{\pi_i}$ contiguously $\forall i$ across secondary storage
    devices.
    \STATE Find $\nabla f(x_0)$ and $P_{\pi_1\pi_1}^{-1}, \ldots, P_{\pi_{m} \pi_{m}}^{-1}.$
    \STATE {\bfseries Main Procedure:}
    \FOR{$k = 0, 1, \ldots,$}
    \STATE Find $\beta_{\pi_1}(x_k), \ldots, \beta_{\pi_m}(x_k).$
    \STATE  $\pi \leftarrow \underset{\pi_i \in \Pi}{\text{argmax}} \hspace{0.2em} \beta_{\pi_i}(x_k).$
    \STATE Load $P_{\pi}$ and $q_{\pi}$ from secondary storage.
    \STATE $\alpha \leftarrow (P_{\pi \pi})^{-1} (q_{\pi} - P_{\pi} x_k).$
    \STATE $x_{k + 1} \leftarrow x_k + \mathbb{I}_{\pi}^t \alpha.$
    \STATE $\nabla f(x_{k + 1}) \leftarrow \nabla f(x_k) + P_{\pi}^t \alpha.$
    \ENDFOR
\end{algorithmic}
\end{algorithm}

We end this section with a brief discussion on why an equivalent greedy strategy for BK method from Algorithm~\ref{alg:BKSch} is not HDC. Let $\{x_k\}$ be the iterates of the BK method. Then, from Lemma~\ref{lem:BKisDLDR} and Fact~\ref{fact:GuaDesc}, it follows that $||x_k - x_{\text{opt}}||_2$ is a non-increasing function of $k,$ whatever be the input partition $\Pi$ and strategy $\zeta_b.$  In fact, if block $P_{\pi_i}$ is chosen in the $k-$th iteration, then
\begin{equation}
\label{eqn:PythThmNorm2}
||x_{k + 1} - x_{\text{opt}}||_2^2 = ||x_k - x_{\text{opt}}||_2^2 - (q_{\pi_i} - P_{\pi_i} x_k)^t (P_{\pi_i} P_{\pi_i}^t)^{-1} (q_{\pi_i} - P_{\pi_i}x_k).
\end{equation}
Based on this equation, the greedy strategy for BK methods would be to pick that block $P_{\pi_i}$ for which $(q_{\pi_i} - P_{\pi_i} x_k)^t (P_{\pi_i} P_{\pi_i}^t)^{-1} (q_{\pi_i} - P_{\pi_i}x_k)$ is largest amongst all possible choices. Here again, it follows that the knowledge of the complete gradient vector would be needed in every iteration. But, as mentioned before, successive iterates in BK methods generically differ in all coordinates. Hence, computing gradient in each iteration of a BK method is essentially an $O(n^2)$ operation. This violates the requirements of feature $\mathcal{F}_2$ thereby proving that the greedy BK method is non-HDC. For this reason, we will not pursue this method further.

\section{Performance Analysis of GBCD Method}
\label{sec:PerMeas}
Consider a high dimensional UQP of the form given in \eqref{eqn:UQP}. To solve this using the proposed GBCD method given in Algorithm~\ref{alg:BCDSchGreedy}, observe that one can use different partitions of $P$ as input. Our goal here is to study the effect of the input partition $\Pi$ on the GBCD method's convergence rate and total running time to find an approximate solution of the given UQP. This is needed to check if the GBCD method possesses feature $\mathcal{F}_4.$

Let $\{x_k\}$ be iterates of the GBCD method. Because of Lemma~\ref{lem:BCDMonCov}, it follows that $||x_k - x_{\text{opt}}||_p$ is a non-increasing function of $k,$ whatever be the input partition $\Pi.$ Hence, it makes sense to use
\begin{equation}
\label{eqn:calC}
\mathcal{C}(\Pi) := \max_{k \geq 0} \mathcal{C}_k(\Pi)
\end{equation}
to define the convergence rate of the GBCD method,
where
\begin{equation}
\label{eqn:calCk}
\mathcal{C}_k(\Pi) := \frac{||x_{k + 1} - x_{\text{opt}}||_P^2}{||x_k - x_{\text{opt}}||_P^2}.
\end{equation}
With $\mathbb{I}_{\pi_i}$ as defined in Algorithm~\ref{alg:BCDSch}, let $\mathbb{I}_{\Pi} \equiv \left[\mathbb{I}_{\pi_1}^t \cdots \; \mathbb{I}_{\pi_{m}}^t\right]^t$ denote a row permutation of the identity matrix $\mathbb{I}.$ Let
\begin{equation}
\label{eqn:PermP}
P_{\Pi} := \mathbb{I}_{\Pi} P \mathbb{I}_{\Pi}^t
\end{equation}
denote the rearrangement of $P$ dictated by $\Pi$ and let
\begin{equation}
\label{eqn:BlkDiag}
B_{\Pi} := \left[
\begin{array}{cccc}
P_{\pi_1\pi_1} & 0 & \cdots & 0\\
0 & P_{\pi_2\pi_2} & 0 & 0\\
\vdots & \vdots  & \ddots & \vdots\\
0 & 0 & \cdots & P_{\pi_{m}\pi_{m}}
\end{array}
\right]
\end{equation}
denote a block diagonal matrix made up of the block diagonal entries of $P_{\Pi}.$ Note that $P_{\Pi} \in \mathbb{S}_{++}^{n}$ as well as $B_{\Pi} \in \mathbb{S}_{++}^{n}.$ We prove two results before obtaining a bound on convergence rate in Lemma~\ref{lem:conRate}.

\begin{lemma}
\label{lem:minEigValPos}
Eigenvalues of $P_{\Pi}B_{\Pi}^{-1}$ are all real and $\lambda_{\min} (P_{\Pi} B_{\Pi}^{-1}) > 0.$
\end{lemma}
\proof{Proof.}
Observe that $P_{\Pi}B_{\Pi}^{-1}$ and $B_{\Pi}^{-1/2} P_{\Pi} B_{\Pi}^{-1/2}$ have the same set of eigenvalues. The desired result now follows since $B_{\Pi}^{-1/2} P_{\Pi} B_{\Pi}^{-1/2} \in \mathbb{S}_{++}^n.$ \Halmos
\endproof

\begin{lemma}
\label{lem:minRatVal}
For $y \in \mathbb{R}^n,$ let $\psi(y) := y^tB_{\Pi}^{-1}y/y^tP_{\Pi}^{-1}y.$ Then for each $y$
\[
\psi(y) \geq \lambda_{\min}(P_{\Pi} B_{\Pi}^{-1}).
\]
\end{lemma}
\proof{Proof.}
Observe that
\[
\psi'(y) = \frac{(y^tP_{\Pi}^{-1}y) B_{\Pi}^{-1}y - (y^tB_{\Pi}^{-1}y) P_{\Pi}^{-1}y}{(y^tP_{\Pi}^{-1}y)^2}.
\]
Setting $\psi'(y) = 0,$ it follows that the extrema of the function $\psi$ occurs at precisely those $y$ at which $
P_{\Pi} B_{\Pi}^{-1} y = \psi(y) y,$ i.e, $y$ is the eigenvector of $P_{\Pi} B_{\Pi}^{-1}.$ Furthermore, the extremum values are the eigenvalues of $P_{\Pi} B_{\Pi}^{-1}.$ The desired result thus follows. \Halmos
\endproof

\begin{lemma}
\label{lem:conRate}
$\mathcal{C}(\Pi) \leq 1 - \frac{1}{m} \lambda_{\min} (P_{\Pi} B_{\Pi}^{-1}).$
\end{lemma}
\proof{Proof.}
Let $k$ be an arbitrary but fixed iteration index. Suppose that
\[
\pi = \underset{\pi_i \in \Pi}{\text{argmax}} \beta_{\pi_i}(x_k).
\]
Then, clearly
\[
\beta_{\pi}(x_k) \geq \frac{1}{m} \sum_{\pi_i \in \Pi} \beta_{\pi_i}(x_k).
\]
Hence from \eqref{eqn:PythThmNrmP}, we have that
\[
\mathcal{C}_k (\Pi) \leq 1 - \frac{1}{m} \frac{\sum_{\pi_i \in \Pi} \beta_{\pi_i}(x_k)}{||x_{k} - x_{\text{opt}}||_P^2}.
\]
Since $\mathbb{I}^t_{\Pi}\mathbb{I}_{\Pi} = \mathbb{I}_{\Pi} \mathbb{I}^t_{\Pi} = \mathbb{I},$ observe that if $y_k := \mathbb{I}_{\Pi}P(x_k - x_{\text{opt}}),$ then $||x_k - x_{\text{opt}}||_P^2 = y_k^t P_{\Pi}^{-1} y_k.$ Further, $\sum_{\pi_i \in \Pi} \beta_{\pi_i}(x_k) = y_k^t B_{\Pi}^{-1} y_k.$ Putting all this together, we get
\[
\mathcal{C}_k(\Pi) \leq 1 - \frac{1}{m}\frac{y_k^t B_{\Pi}^{-1} y_k}{y_k^t P_{\Pi}^{-1} y_k}.
\]
The desired result now follows from Lemma~\ref{lem:minRatVal}. \Halmos
\endproof
From Lemmas~\ref{lem:minEigValPos} and \ref{lem:conRate}, we have the following result.
\begin{corollary}
\label{cor:convRateLessThan1}
 $\mathcal{C}(\Pi) < 1.$
\end{corollary}

The last two results prove the following theorem.

\begin{theorem}
\label{thm:convGBCD}
Let $\{x_k\}$ be the iterates of the GBCD method given in Algorithm~\ref{alg:BCDSchGreedy}. Then, $x_k \rightarrow x_{\text{opt}}$ as $k \rightarrow \infty$ with a convergence rate bounded above by $1 - \frac{1}{m} \lambda_{\min}(P_{\Pi}B^{-1}_{\Pi}).$
\end{theorem}

We now use this result to obtain a bound on the total running time (number of iterations $\times$ running time per iteration) of the GBCD method to find an approximate solution of the given UQP. Fix $\epsilon > 0.$ For the given initial approximation $x_0$ of $x_{\text{opt}}$ in GBCD method, we will say that the $k-$th estimate $x_k$ is $\epsilon-$close to $x_{\text{opt}}$ if
\begin{equation}
\label{eqn:epsCloseSoln}
\mathcal{E}_k := \frac{||x_k - x_{\text{opt}}||_P}{||x_0 - x_{\text{opt}}||_P } < \epsilon.
\end{equation}
Clearly, $\mathcal{E}_k \leq \mathcal{C}(\Pi)^{k/2}.$ From this, it follows that if $k \geq 2\log(\frac{1}{\epsilon})/\log(\mathcal{C}(\Pi)),$ then $x_k$ is certainly an $\epsilon-$close solution. Hence, from Lemma~\ref{lem:conRate} and the fact that the per iteration running time is $O(nd),$ it follows that the total running time of the GBCD method to find an $\epsilon-$close solution is $O(n d m\log(1/\epsilon)/\lambda_{\min}(P_{\Pi} B^{-1}_{\Pi})).$ Now note that all eigenvalues of $B_{\Pi}$ lie between $\lambda_{\min}(P)$ and $\lambda_{\max}(P)$ and eigenvalues of $P_{\Pi}$ are same as those of $P.$ Hence, it follows that, for any partition $\Pi$, $1/\lambda_{\min}(P_{\Pi} B_{\Pi}^{-1}) \leq \kappa(P).$ Recall that for the steepest descent method the total running time to find an $\epsilon-$close solution is $O(n^2\kappa(P) \log(1/\epsilon)).$ From these, it follows that the total running time of the GBCD method is comparable to that of the steepest descent. Simulations in fact show that the GBCD is typically much better off than the steepest descent and the other BCD methods. Loosely speaking, this shows that the GBCD method does indeed possess feature $\mathcal{F}_4.$

\section{Good Partitioning Strategies}
\label{sec:ParStra}
Consider a high dimensional version of the UQP in \eqref{eqn:UQP}. Theorem~\ref{thm:convGBCD} tells that when the GBCD method of Algorithm~\ref{alg:BCDSchGreedy} is used to solve this UQP, its convergence rate is influenced by the partition of $P$ that is given as input. This suggests that a good partition $\Pi \in \wp$ may ensure faster convergence for the GBCD method and hence lesser time to find an approximate solution of the given UQP. We briefly discuss here what constitutes a good partition and suggest heuristic ways to find it.

Fix $\epsilon >  0.$ Let $T_{\epsilon}(\Pi)$ denote the time taken by GBCD method, with partition $\Pi$ as input, to find an $\epsilon-close$ solution of given UQP. We will say that a row partition $\Pi^* \in \wp$ of $P$ is good if:

\begin{enumerate}
\item $T_{\epsilon} (\Pi^*)$ is less than or close to $T_{\epsilon} \left(\underset{\Pi \in \wp_d}{\text{argmin}} \hspace{1ex} C(\Pi) \right),$ where $C(\Pi)$ is as in \eqref{eqn:calC},

\item there exists a HDC method to find $\Pi^*$ and

\item the time taken to find $\Pi^*$ is small relative to the time taken by the GBCD method to solve the given UQP using an arbitrary partition of $P$ as input.
\end{enumerate}

Finding such a good partition $\Pi^*$ at present seems hard and we leave it as a future objective. What we do next, instead, is to come up with a heuristic way to find partitions that ensure better convergence for the GBCD method. In this direction, we first simplify the bound on the convergence rate. For a matrix $A,$ let $||A||_2$ denote its spectral norm.

\begin{lemma}
\label{lem:betConv}
Let $P_{\Pi}$ and $B_{\Pi}$ be as in \eqref{eqn:PermP} and \eqref{eqn:BlkDiag}. Then,
\[
\lambda_{\min} (P_{\Pi} B_{\Pi}^{-1}) \geq 1 - ||P_{\Pi} - B_{\Pi}||_2/\lambda_{\min}(B_{\Pi}).
\]
\end{lemma}
\proof{Proof.}
Observe that
\begin{eqnarray*}
||\mathbb{I} - P_{\Pi}B_{\Pi}^{-1}||_2 & = & ||(B_{\Pi} - P_{\Pi}) B_{\Pi}^{-1}||_2 \\
& \leq & \frac{||B_{\Pi} - P_{\Pi}||_2}{\lambda_{\min}(B_{\Pi})}.
\end{eqnarray*}
This implies that for each eigenvalue $\lambda$ of $P_{\Pi} B_{\Pi}^{-1}$
\[
|1 - \lambda| \leq \frac{||B_{\Pi} - P_{\Pi}||_2}{\lambda_{\min}(B_{\Pi})}.
\]
The desired result now follows. \Halmos
\endproof

\begin{corollary}
\label{cor:heuStr}
Let $P_{\Pi}$ and $B_{\Pi}$ be as in \eqref{eqn:PermP} and \eqref{eqn:BlkDiag} and $C(\Pi)$ as in \eqref{eqn:calC}. Then
\[
C(\Pi) < 1 - \frac{1}{m}\left(1 - \frac{||P_{\Pi} - B_{\Pi}||_2}{\lambda_{\min}(B_{\Pi})}\right).
\]
\end{corollary}
\proof{Proof.}
This follows from Lemmas~\ref{lem:conRate} and \ref{lem:betConv}. \Halmos
\endproof

Although weaker than Lemma~\ref{lem:conRate}, the bound given in Corollary~\ref{cor:heuStr} gives a simpler understanding of the influence of $\Pi$ on the convergence rate. Loosely speaking, the new bound says that closer $B_\Pi$ is to $P_\Pi,$ faster is the convergence rate. In fact, note that if $P_\Pi = B_{\Pi},$ then the convergence rate is bounded above by $1 - 1/m;$ a number of independent of the eigenvalues of $P.$ Based on these observations, we suggest the following heuristic idea to improve the convergence rate of the GBCD method: \emph{Pick a partition $\Pi \in \wp$ for which $||P_{\Pi} - B_{\Pi}||_2$ is small and $\lambda_{\min}(B_{\Pi})$ is large.} We will call such a partition of $P$ as block diagonally dominant. Experiment 3 of Section~\ref{sec:SimRes} gives an example where using a block diagonally dominant partition indeed speeds up convergence.

\section{Parallel and Distributed Implementation}
\label{sec:ParDisImp}
By a parallel computing system, we imply a system which has multiple processors that share a common main memory and a set of secondary storage devices. Distributed systems, on the other hand, will mean a group of networked computers each equipped with a single processor, main memory and secondary storage device of its own. Our objective here is to discuss briefly the necessary changes to be made in the GBCD method of Algorithm~\ref{alg:BCDSchGreedy} to take advantage of the additional resources available in these two setups. For the UQP to be solved, we will, as usual, use $d$ to denote the maximum of the number of rows in the chosen partition of $P.$ For pedagogical considerations, we will assume that $x_0 = 0$ throughout this section. As mentioned before, the advantage is that the initial gradient $\nabla f(x_0)$ is readily available.

\subsection{Parallel Implementation}
Let $1 < N_p \leq d$ be the number of available parallel processors. Observe that the preprocessing phase involves computing the inverse of $m$ matrices each of dimension at most $d \times d.$ The first modification we then suggest is to divide these set of matrices into $N_p$ groups of roughly equal sizes and feed one group as input to each processor. Then compute the inverse of matrices in parallel. Clearly, the time required to compute all the inverses will reduce from $O(nd^2)$ to $O(nd^2/N_p).$

Next observe that each iteration of the main  procedure involves three phases: first is computing $\pi,$ second is reading $P_{\pi}$ and $q_{\pi}$ from secondary storage into main memory, and third is determining $x_{k + 1}$ from $x_k.$ Based on these, we suggest two modifications. Since the inverse of $P_{\pi \pi},$ $\forall \; \pi \in \Pi,$ is already known, observe that the first and the third phase only involve matrix vector multiplications. Carry out these operations by dividing the rows of the concerned matrix across the processors while giving the vector involved to every processor. For phase two, read the entries of $P_{\pi}$ and $q_{\pi}$ from secondary storage into main memory through $N_p$ streams in parallel. Note that this latter idea suggesting parallelization in data fetch from secondary storage may perhaps require new hardware technology and it is one of our future objectives to understand this operation in more detail. With these modifications, the run time per iteration will reduce from $O(nd)$ to $O(nd/N_p).$ Consequently, the total running time to find an $\epsilon-$close solution will reduce to $O(nmd \log(1/\epsilon)/(N_p \lambda_{\min} (P_{\Pi} B_{\Pi}^{-1}))).$

\subsection{Distributed Implementation}
In the distributed model, let us suppose that there are $N_p$ independent computer nodes connected with each other. In this setup, the modification to Algorithm~\ref{alg:BCDSchGreedy} we suggest is the following. First divide the $m$ blocks $P_{\pi_1}, \ldots, P_{\pi_m}$ into $N_p$ groups and associate each group to one of the nodes. Store the blocks in the associated node's secondary storage device. Store the matrices $P_{\pi_1 \pi_1}^{-1}, \ldots, P_{\pi_m \pi_m}^{-1},$ however, at every node's secondary storage device. In fact, these can be retained in the main memory associated with each node throughout the running time of the algorithm. Recall that these matrices only need $O(nd)$ space. These modifications concern the preprocessing phase.

The main procedure is to be implemented as follows. Input $x_0$ and $\nabla f(x_0)$ to an arbitrary node and compute $\pi$ as defined in Algorithm~\ref{alg:BCDSchGreedy}. Now inductively, for $k = 0, 1, \ldots,$ do the following.
\begin{enumerate}
\item Transfer $x_k$ and $\nabla f(x_k)$ to the node which holds $P_{\pi}$ and $q_{\pi}.$

\item At this node, compute $x_{k + 1}$ and $\nabla f(x_{k + 1})$ and also determine the new $\pi.$
\end{enumerate}

In the first step, note that we only need to transfer $2n$ entries. Hence, the time required for this operation will be $O(n).$ The second step, as usual, will require $O(nd)$ time. We emphasize here that the above implementation is serial. That is, at any given time, precisely one node is active while the rest of the nodes do no computation. This may lead to the belief that this idea has very poor resource utilization. As we  show in the two scenarios below, this need not be the case always. In fact, the above distributed implementation scheme may be very advantageous. The first scenario is when the cost of the delay incurred in transmitting $P_{\pi}$ and $q_{\pi}$ from the secondary storage device to main memory far outweighs the cost of having $N_p$ independent computing devices (which we call here as nodes). This may happen, for example, if the secondary storage devices are spread out in different geographic locations and we have a centrally located processor. The second scenario is when multiple high dimensional UQPs defined using the same $P$ matrix but different $q$ vectors need to be solved. By assigning a copy of GBCD method for each of the UQPs, one can solve all of them simultaneously within the same distributed setup. If these copies concurrently work on different nodes in most of the iterations, then efficient utilization of resources can be achieved. We do not pursue this here.

\section{Simulation Results}
\label{sec:SimRes}
In this section, we give details of the results obtained from three simulation experiments done using Matlab. Through this, we wish to demonstrate the utility of the proposed GBCD method from Algorithm~\ref{alg:BCDSchGreedy} in comparison to existing approaches. Specifically, the first experiment gives an example of settings where the proposed GBCD method is a better choice than the greedy BCD of \cite{bo_uai08}, randomized BCD of \cite{nesterov2012efficiency}, the steepest descent and, in fact, even the conjugate gradient method. This experiment also gives a realistic understanding of the difficulties involved in solving high dimensional UQPs. The second experiment highlights the limitations of the static block selection strategies used by methods such as the randomized coordinate descent of \cite{leventhal2010randomized} and randomized Kaczmarz from \cite{strohmer2009randomized}. Specifically, it talks of a scenario where convergence of these randomized algorithms is good only with respect to error in $||\cdot||_P$ norm or equivalently function values (see \eqref{eqn:fValNormPRel}) but not in $||\cdot||_2$ norm. Intuitive reasons for why the proposed GBCD method will not have such behaviour is also given here. The third experiment exhibits a scenario where using a block diagonally dominant partition, as was defined in Section~\ref{sec:ParStra}, is better than using an arbitrary partition. Note that, in the second and third experiment, our focus is more on getting the idea across than on implementation issues of high dimensional UQPs. Hence we work with only moderately sized input matrices. The technical specifications of the machine on which the experiments were carried out are as follows. The machine had $3$ GB RAM, $60$ GB secondary storage space and an Intel i3 processor with $2.6$ GHz clock cycle.

\subsection{Experiment 1}
Here, we chose $n = 2^{15}= 32768.$ To generate the matrix $P \in \mathbb{S}_{++}^{n},$ we used the following logic. We first generated a random matrix $V \in \mathbb{R}^{n \times n}$ and then used the relation $P = V^tV$ to generate $P.$ But note that the size of both $P$ and $V$ equals $8$GB here.  Since only $3$GB RAM was available, we had to generate both these matrices by breaking them into smaller submatrices and individually generating these submatrices. Formally, if $(V_{ij})_{1 \leq i, j \leq 256},$ where each $V_{ij} \in \mathbb{R}^{128 \times 128},$ denotes a block partition form of $V,$ then we first generated these $256^2$ submatrices $\{V_{ij}\}$ using the rule
\[
V_{ij} =
 \begin{cases}
 10 \; \mathcal{Z}_{128} & \text{if } i = j \\
 0.1 \; \mathcal{Z}_{128} & \text{otherwise}
 \end{cases},
\]
where $\mathcal{Z}_{128}$ denotes a $128 \times 128$ matrix made up IID standard Gaussian random variables. This rule ensured that the submatrices $\{V_{ii} : 1 \leq i \leq 256\}$ closer to diagonal of $V$ had numerically larger mass relative to the other submatrices. Each of these $256^2$ submatrices were stored separately in individual files on hard disk. To generate these submatrices, we required approximately $11$ minutes. Let $P_{\pi_1} \in \mathbb{R}^{128 \times 32768}$ denote the first $128$ rows of $P,$ $P_{\pi_2}$ the next $128$ rows of $P$ and so on. Using the submatrices of $V,$ we then generated the blocks $P_{\pi_1}$ to $P_{\pi_{256}}$ and stored each of them individually in a separate file. This entire operation to generate the $P$ matrix took approximately $4$ hours and resulted in $256$ files each of size $32$ MB. Note that entries in the submatrices $\{P_{\pi_i\pi_i} : 1\leq i \leq 256\},$ where $P_{\pi_i \pi_i}$ is as defined in Algorithm~\ref{alg:BCDSch}, are larger numerically compared to other entries of $P.$ In other words, if $\Pi \equiv \{\pi_1, \ldots, \pi_{256}\},$ where $\pi_1 = \{1,\ldots, 128\},$ $\pi_2 = \{129, \ldots, 256\}$ and so on, then $\Pi$ denotes a block diagonally dominant partition of $P$ in the sense described in Section~\ref{sec:ParStra}.  We also computed the inverses $\{P_{\pi_i\pi_i}^{-1}: 1 \leq i \leq 256\}$ and stored these separately in another $32$ MB file. This step took roughly $3.5$ minutes. An arbitrary optimal point $x_{\text{opt}}$ was chosen and the vector $q$ was computed using the relation $q = Px_{\text{opt}}.$

\begin{figure}[t!]
\begin{center}
\centerline{\includegraphics[scale = 0.75]{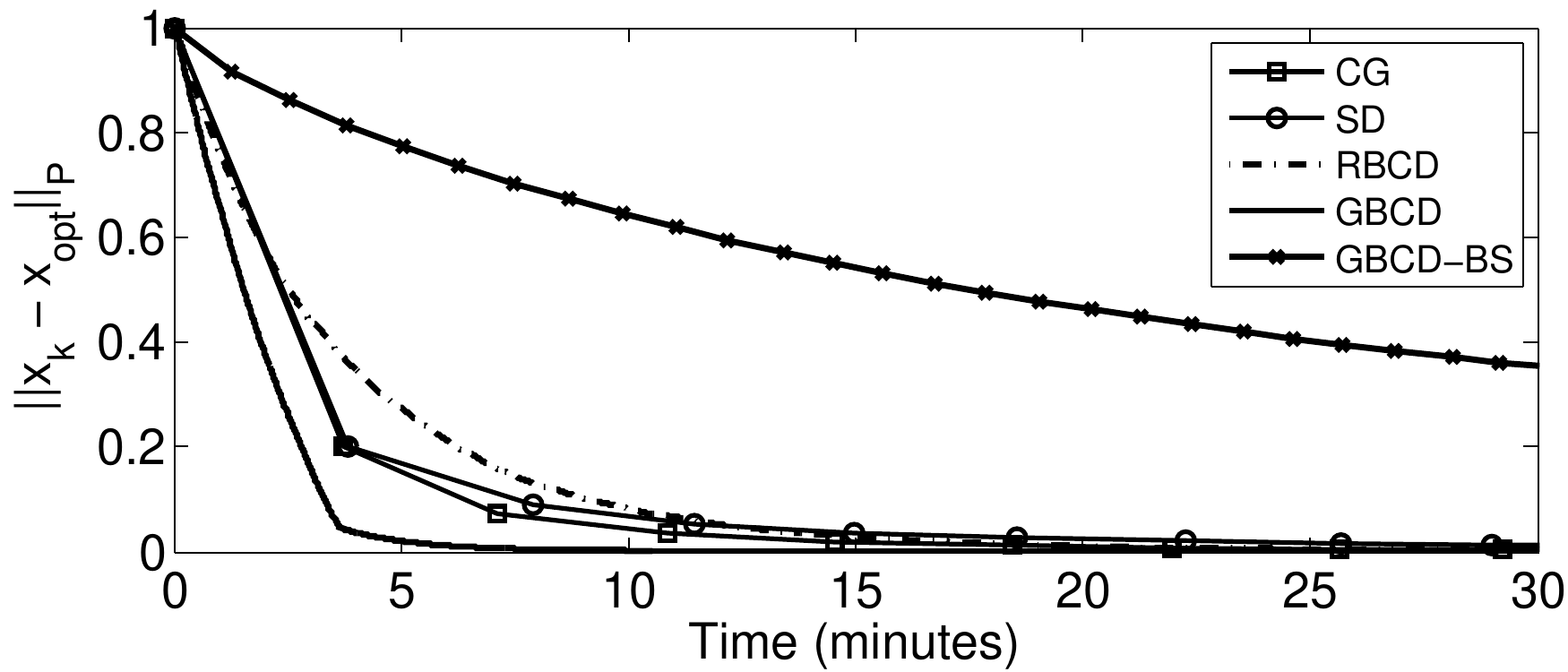}}
\caption{Comparison of standard algorithm for almost block diagonal partition.}
\label{fig:Expt1}
\end{center}
\end{figure}

We solved a UQP with $P$ and $q$ as given in the above setup using five methods: the conjugate gradient, steepest descent, Nesterov's randomized BCD, Bo and Sminchisescu's greedy BCD and the proposed GBCD of Algorithm~\ref{alg:BCDSchGreedy}. To implement Algorithm~\ref{alg:BCDSchGreedy}, the partition $\Pi$ defined earlier was given as input to the GBCD method. The randomized BCD method was implemented exactly as Algorithm~\ref{alg:BCDSchGreedy} except for the block selection strategy. Specifically, the blocks of the partition were randomly chosen using a distribution proportional to the maximum eigenvalues of the $128 \times 128$ matrices $P_{\pi_1\pi_1}, \ldots, P_{\pi_{m} \pi_{m}}.$ The greedy BCD of \cite{bo_uai08} was implemented using Algorithms 1 and 2 given in that paper. This method chose $128$ rows in each iteration. For the conjugate gradient and steepest descent methods, recall that each iteration requires computation of the gradient. We did each such gradient computation in $256$ stages where the $i-$th stage dealt with reading the row submatrix $P_{\pi_i}$ and computing the partial gradient (see \eqref{eqn:partGrad}) associated with it. The initial estimate for all these methods was the origin. The error under $||\cdot||_P$ norm (scaled to start at 1) versus time plot for the four methods is given in Figure~\ref{fig:Expt1} with CG, SD, RBCD, GBCD-BS and GBCD identifying the individual methods. Note that the output of the randomized BCD method was averaged over $25$ runs. As one can see, the proposed GBCD method converges faster than the other methods. This demonstrates the intended superiority of the proposed approach in settings where the input matrix is block diagonally dominant.

The plot in Figure~\ref{fig:Expt1} also highlights the crucial feature of low running times per iteration for the proposed GBCD and the randomized BCD methods. We have placed the marker {\tiny $\mathbf{\Box}$}, $\circ$ and {\tiny{$\times$}} on the progress trajectories of conjugate gradient, steepest descent and GBCD-BS precisely at times where we obtained their respective estimates. As one can see, each iteration of the conjugate gradient and steepest descent method, on an average, took roughly $3$ minutes $45$ seconds. The average iteration for the GBCD-BS method took roughly $70$ seconds. In sharp contrast, each iteration of the GBCD and RBCD methods took only $0.8$ seconds. In other words, during the $30$ minute period that all these algorithm were run, the steepest descent and conjugate gradient methods gave $9$ estimates, the GBCD-BS gave $26$ estimates while the GBCD and RBCD methods resulted in over $2000$ estimates.

\subsection{Experiment 2}
For this experiment, we chose $n = 1024.$ We first generated an arbitrary matrix $V \in \mathbb{R}^{n \times n}$ using the rule $V = \mathcal{Z}_{1024}$ and then constructed the matrix $\tilde{P}$ using the relation $\tilde{P} = V^t V.$  We then chose a subset $\tau \subset [n]$ made of $32$ arbitrary indices. The rows and columns of $\tilde{P}$ corresponding to $\tau$ were multiplied by $1000$ to finally get the $P$ matrix. Note that, since $n$ was only $1024,$ we stored all the above matrices in their entirety directly in the main memory. Finally we generated an arbitrary $x_{\text{opt}}$ and built $q$ using the relation $q = P x_{\text{opt}}.$

\begin{figure}[t!]
\begin{center}
\centerline{\includegraphics[scale = 0.75]{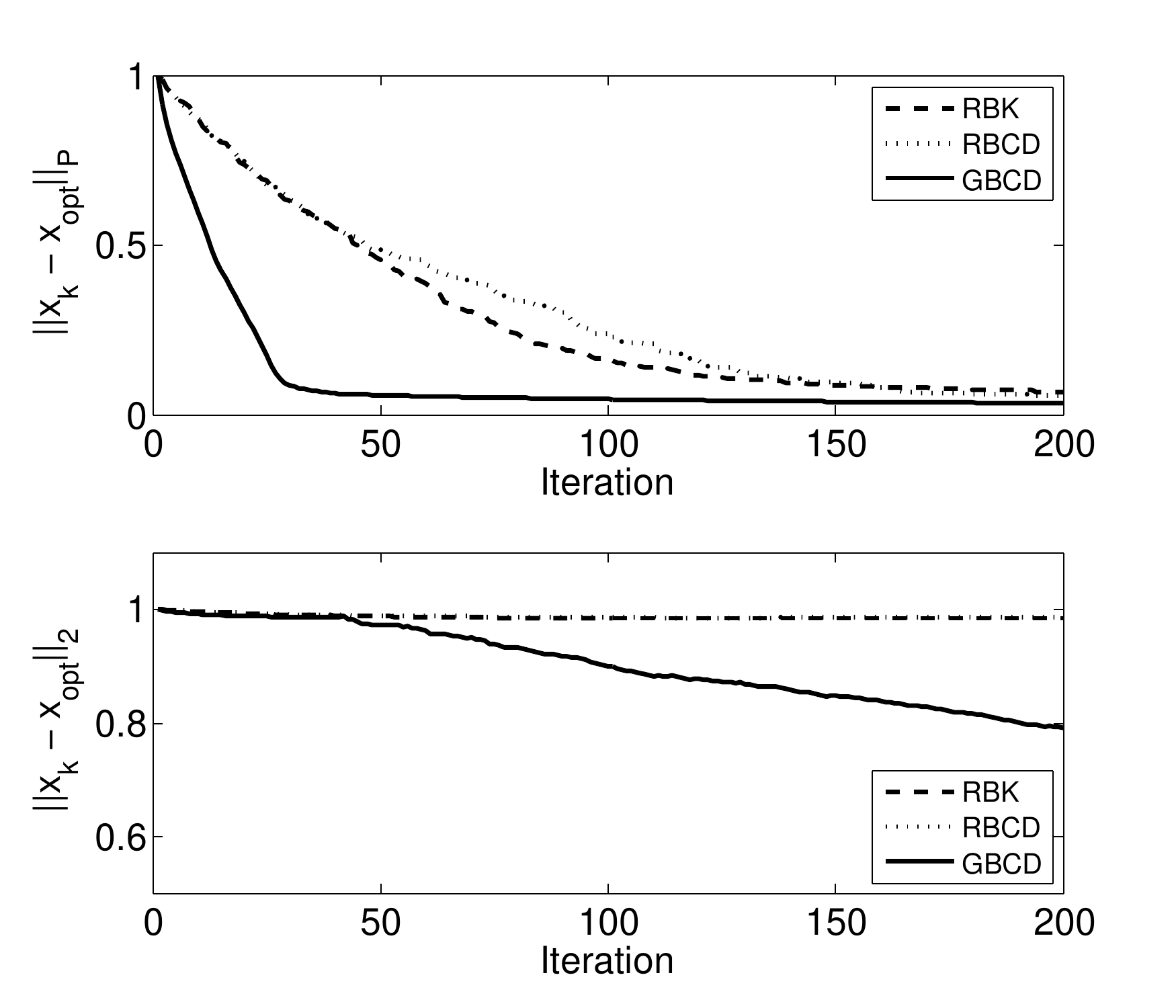}}
\caption{Randomized versus greedy strategy.}
\label{fig:Expt2}
\end{center}
\end{figure}

We solved a UQP with this $P$ and $q$ using three methods: Strohmer and Vershynin's randomized Kaczmarz, Leventhal and Lewis' randomized coordinate descent and the proposed GBCD. We implemented the GBCD method using Algorithm~\ref{alg:BCDSchGreedy} with $\Pi \equiv \{\pi_1, \ldots, \pi_{1024}\}$ where each $\pi_i = \{i\}.$ That is, each $P_{\pi_i}$ is made of one row of $P.$ Since the input matrix $P$ was small enough to be retained in the main memory, we skipped all the steps of this algorithm that involved secondary storage reads. The computation of the inverses in the preprocessing step was also skipped. The randomized coordinate descent was implemented in exactly the same manner as the GBCD method except for the block selection strategy. In each iteration of the randomized BCD, the blocks from $\Pi$ (rows of $P$ in this case) were chosen using a distribution that was proportional to the diagonal entries of $P.$ Along similar lines, we implemented the randomized Kaczmarz method using Algorithm~\ref{alg:BKSch}. Specifically, the strategy $\zeta_b$ used was to select the blocks of the partition $\Pi$ using a distribution that was proportional to the square of $||\cdot||_2$ norm of the rows. All the methods started out at the origin. The comparative performance of these methods is given in Figure~\ref{fig:Expt2}. The three methods are denoted using RK, RCD and GBCD. Clearly, under both the norms $||\cdot||_2$ and $||\cdot||_P,$ the GBCD method has faster convergence. But the key thing to observe is that the decrease in error under $||\cdot||_2$ norm for the randomized approaches is almost negligible. This happens mainly because of the way the $P$ matrix was constructed and the way the randomized methods work. In particular, note that both these methods sample almost always only those rows of $P$ that have large mass, i.e., the rows of the row submatrix $P_{\tau},$ irrespective of where the current estimate is. Because of this, the estimates of the Kaczmarz algorithm lie very close to the subspace $\text{span}{(P_{\tau})}.$ In similar fashion, the estimates of the coordinate descent algorithm lie very close to the subspace $\text{span}{(\mathbb{I}_{\tau})}.$ From these observation, it is easy why the $||\cdot||_2$ norm convergence of the error for these methods is poor. The proposed method, however, is adaptive and hence is able to overcome these problems.

\subsection{Experiment 3}
\label{subsec:Expt3}
\begin{figure}[t!]

\begin{center}
\centerline{\includegraphics[scale = 0.75]{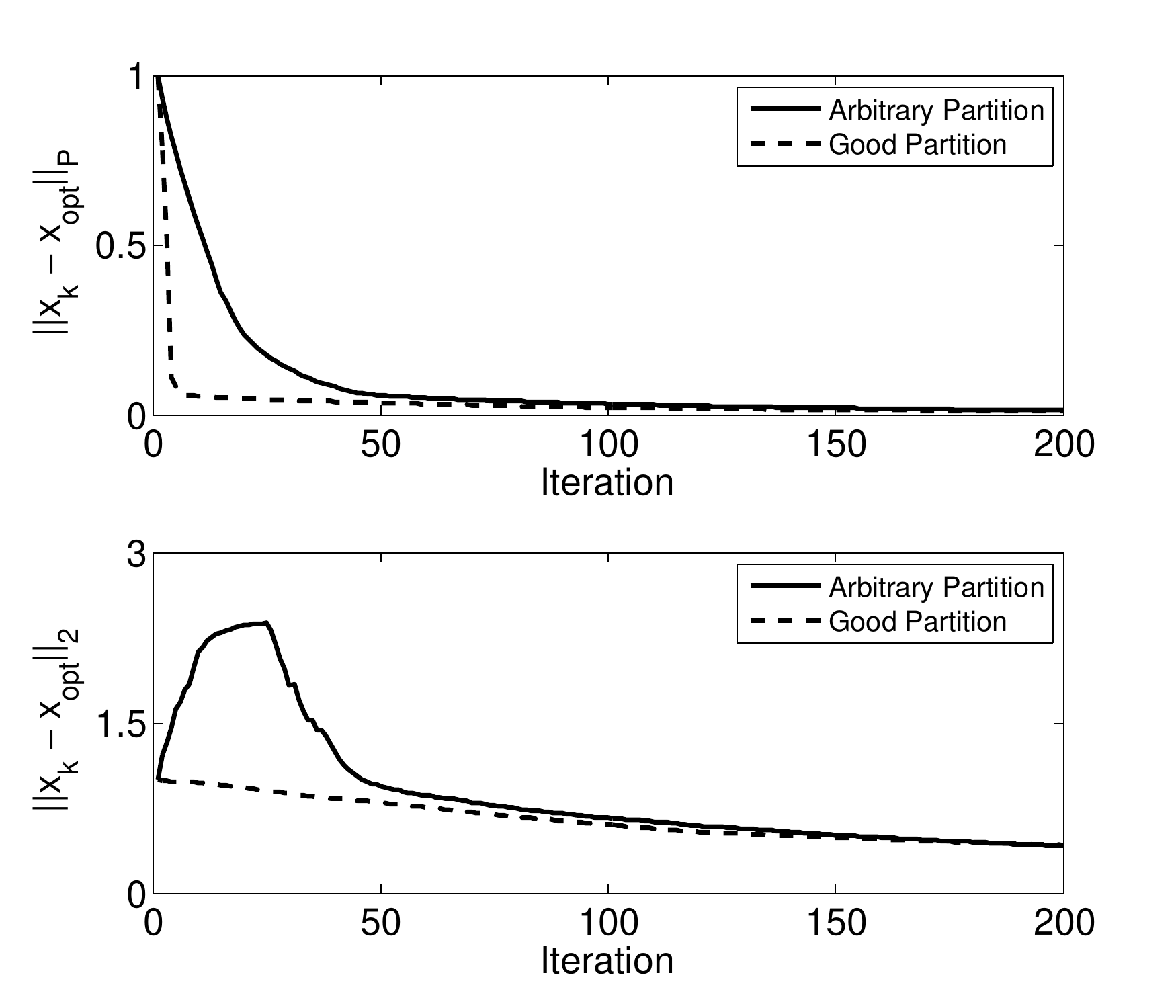}}
\caption{Choosing a good partition.}
\label{fig:Expt3}
\end{center}
\end{figure}

In this experiment, we chose $n = 1024$ and generated the $P \in \mathbb{R}^{n \times n}$ matrix and $q \in \mathbb{R}^n$ vector in exactly the same manner as we did in Experiment $2.$ We then solved a UQP with this $P$ matrix and $q$ vector as input using the proposed GBCD method in two different ways. Specifically, we chose two different input partitions. In the first way, we arbitrarily partitioned the rows of $P$ into $32$ blocks with each block made up of $32$ rows. In the second way, we kept all the $32$ rows with large numbers, i.e., rows with indices in $\tau,$ in one block while the remaining rows were arbitrarily partitioned into $31$ blocks with each block again made up of $32$ rows. In a loose sense, we tried to use a block diagonally dominant partition of $P.$ The starting point in both ways was the origin. The performance of GBCD with these two different partitions as input is given in Figure~\ref{fig:Expt3}. As one can see, the performance was better with the block diagonally dominant partition.

\section{Conclusion}
\label{sec:Concl}
In this paper, we discussed solution methods for high dimensional UQPs when one has access only to limited computational resources. We began with desired features of an algorithm in such settings. It was shown that the traditional BK and BCD methods using a prefixed row partition of the input matrix  are the only methods amongst existing ones that can be readily made to possess these features. The key contribution of this paper is the GBCD method from Algorithm~\ref{alg:BCDSchGreedy}. Theoretical and experimental analysis of its convergence rate revealed that using a block diagonally dominant row partition of the input matrix speeds up convergence. In fact, for input matrices which are almost block diagonal, simulation results showed that the proposed method converged faster than all existing methods. Finding the partition under which the convergence will be the fastest, however, remains an open question.

\ACKNOWLEDGMENT{The research of G. Thoppe was supported in part by an IBM Fellowship. The research of V. Borkar was supported in part by an IBM SUR Award, a J.~C.~Bose Fellowship, and a grant for `Distributed Computation for Optimization over Large Networks and High Dimensional Data Analysis' from the Department of Science and Technology, Government of India. A portion of this work was done when G. Thoppe did a summer internship at IBM with D. Garg.}

\bibliographystyle{ormsv080} 
\bibliography{OR_Bib} 

\end{document}